\documentclass{amsart}

\begin{document}

\newcommand{\nc}{\newcommand}

\nc{\pr}{\noindent{\em Proof. }}
\nc{\g}{\mathfrak g}
\renewcommand{\k}{\mathfrak k}
\nc{\A}{\mathcal A}
\nc{\F}{\mathcal F}
\renewcommand{\H}{\mathfrak H}

\newtheorem{theorem}{Theorem}{}
\newtheorem{lemma}[theorem]{Lemma}{}
\newtheorem{corollary}[theorem]{Corollary}{}
\newtheorem{conjecture}[theorem]{Conjecture}{}
\newtheorem{proposition}[theorem]{Proposition}{}
\newtheorem{axiom}{Axiom}{}
\newtheorem{remark}[theorem]{Remark}{}
\newtheorem{example}{Example}{}
\newtheorem{exercise}{Exercise}{}
\newtheorem{definition}{Definition}{}

\renewcommand{\theremark}{}

\title{An analogue of the operator curl for nonabelian gauge groups and scattering theory}

\author[A. Sevostyanov]{A. Sevostyanov}
\address{ Department of Mathematical Sciences \\
University of Aberdeen  \\ Aberdeen AB24 3UE, United Kingdom}
\email{seva@maths.abdn.ac.uk}


\begin{abstract}
We introduce a new perturbation for the  operator curl related to connections with nonabelian
gauge groups over ${\mathbb{R}}^3$. We also prove that the perturbed operator is unitary equivalent to the operator
curl if the corresponding connection is close enough to the trivial one with respect to a certain topology on the space of connections.
\end{abstract}

\maketitle

\pagestyle{myheadings}

\markboth{A. SEVOSTYANOV}{OPERATOR CURL FOR NONABELIAN GAUGE GROUPS}


\section*{Introduction}

It is well known that the operator $\rm curl$ appears in various fields of mathematics and physics.
Although its spectral theory was developed quite recently (see \cite{Fi,P} and references there).
Interesting spectral problems for the operator $\rm curl$ appear when it is defined in nontrivial
regions of the Euclidean space ${\mathbb{R}}^3$.

As usual, another interesting class of spectral problems for the operator $\rm curl$ would appear
if we
were be able to add a perturbation to it in such a way that the perturbed operator is close,
in some sense, to the original one. A natural operator ${\rm curl}_A$ with potential $A$ that
generalizes the
operator curl in differential geometry is the composition of the Hodge star operator and the
covariant derivative associated to a connection $A$ with a nonabelian gauge group and
defined over ${\mathbb{R}}^3$ (see formula (\ref{defrot}) below). But from the point of view of
spectral theory that operator is quite different from $\rm curl$. For instance,
the scattering theory
for the pair $({\rm curl}_A,{\rm curl}_0)$ can not be properly developed since both operators are
not elliptic.

In this paper we introduce another operator $X_A$ related to nonabelian gauge groups. The operator
$X_A$ is not a differential operator (see formula (\ref{defxa})). In fact, to define it one should
also consider an analogue of the divergence operator $\rm div$ for noncommutative gauge groups (see
formula
(\ref{defdiv})). The operator $\rm div$ is naturally related to the operator $\rm curl$ since
the kernel of $\rm div$ is the
natural domain where the operator $\rm curl$ becomes elliptic.

The spectral properties of the operator $X_A$ introduced in this paper are similar to those of
$\rm curl$. In particular, we prove that if the connection $A$ is small in some sense then the operator
$X_A$ is unitary equivalent to $X_0$. Technically the solution to this problem is achieved by
extending $X_A$ to an elliptic operator that acts in a bigger space. This trick is similar to
that used in \cite{BS1,BS2} for the Maxwell operator.

\vskip 0.3cm
\noindent
{\bf Acknowledgments.}
The author thanks Prof. M. Sh. Birman for useful advises.
I would like also to express my gratitude to  N. Filonov and A. Pushnitski.
The analytic details of this paper would never be completed without their support.



\section{Setup}\label{YMsetup}

Let $K$ be a compact semisimple Lie group, $\k$ its Lie algebra and $\g$ the
complexification of $\k$. We denote by $(\cdot ,\cdot )$ the Killing form of $\g$.
Recall that the restriction of this form to $\k$ is nondegenerate and negatively defined.

Let $\Omega^*(\mathbb{R}^{3},\k)$ be the space of smooth
$\k$-valued differential forms on $\mathbb{R}^{3}$, and $\Omega^*_c(\mathbb{R}^{3},\k)$ the space
of smooth $\k$-valued differential forms on $\mathbb{R}^{3}$ with compact support. We
define a scalar product on  $\Omega^*_c(\mathbb{R}^{3},\k)$ by
\begin{equation}\label{prod}
<\omega_1,\omega_2>=-\int_{\mathbb{R}^{3}}(\omega_1\wedge,*\omega_2)=
-\int_{\mathbb{R}^{3}}*(\omega_1\wedge,*\omega_2)d^{3}x,~
\omega_{1,2} \in \Omega^*_c(\mathbb{R}^{3},\k)
\end{equation}
where $*$ stands for the Hodge star operation associated to the standard
Euclidean metric on  $\mathbb{R}^{3}$, and we evaluate the Killing form on the
values of $\omega_1$ and $*\omega_2$ and also take their exterior product.

Let $A\in \Omega^1(\mathbb{R}^{3},\k)$ be $\k$-valued connection
1-form in the trivial $K$-bundle, associated to the adjoint
representation of $K$, over $\mathbb{R}^{3}$, $F\in
\Omega^1(\mathbb{R}^{2},\k)$ its curvature 2-form, $F=dA +
\frac{1}{2} [A\wedge A]$. Here as usual we denote by $[A\wedge A]$
the operation which takes the exterior product of $\k$-valued
1-forms and the commutator of their values in $\k$, and $d$ stands for the exterior derivative.

We recall that the covariant derivative
$d_A:\Omega^n_c(\mathbb{R}^{3},\k)\rightarrow \Omega^{n+1}_c(\mathbb{R}^{3},\k)$ associated to $A$
is defined by $d_A\omega =d\omega + [A\wedge \omega]$, and the operator
formally adjoint to $d_A$ with respect to scalar product (\ref{prod}) is
equal to $-*d_A*$. We denote by ${\rm div}_A$ the part of this operator acting
from $\Omega^1_c(\mathbb{R}^{3},\k)$ to $\Omega^{0}_c(\mathbb{R}^{3},\k)$, with
the opposite sign,
\begin{equation}\label{defdiv}
{\rm div}_A=*d_A*:\Omega^1_c(\mathbb{R}^{3},\k)\rightarrow
\Omega^{0}_c(\mathbb{R}^{3},\k).
\end{equation}

Consider the affine space  of smooth connections in the trivial
$K$-bundle, associated to the adjoint representation of $K$, over
${\mathbb R}^{3}$. We fix the standard trivialization of this
bundle and the trivial connection as an origin in the affine space
of connections and identify this space with the space
$\Omega^1(\mathbb{R}^{3},\k)$ of $\k$-valued 1-forms on ${\mathbb
R}^{3}$. We shall frequently write $\mathcal{D}$
instead of $\Omega^1(\mathbb{R}^{3},\k)$.

Let $\mathcal{K}$ be the group
of $K$-valued smooth maps $g:\mathbb{R}^{3}\rightarrow K$. $\mathcal{K}$ is called the gauge
group. The Lie algebra of $\mathcal{K}$
is isomorphic to $\Omega^0(\mathbb{R}^{3},\k)$.

The gauge group $\mathcal{K}$ acts on the space of connections $\mathcal{D}$
by
\begin{equation}\label{Gaugeact}
\begin{array}{l}
\mathcal{K}\times \mathcal{D} \rightarrow \mathcal{D},\\
\\
g\times A \mapsto g\circ A= -dgg^{-1}+gAg^{-1},
\end{array}
\end{equation}
where we denote $dgg^{-1}=g^*\theta_R$, $gAg^{-1}=\mathrm{Ad}g(A)$,
$\theta_R$ is the right-invariant Maurer-Cartan form on $K$, and $\mathrm{Ad}g$ stands for $g$ acting in the
adjoint representation.

Action (\ref{Gaugeact}) induces the corresponding transformation law for the curvature $F$,
$$
g\circ F=\mathrm{Ad}g(F),
$$
where $\mathrm{Ad}g$ acts on $F$ componentwise.

For any real vector space $V$ we denote by $V_{\mathbb{C}}$ its
complexification. Let $\mathcal{H}^1$ be the completion of the space $\Omega^1_c(\mathbb{R}^3,\g)$
with respect to
scalar product (\ref{prod}) extended from $\Omega^1_c(\mathbb{R}^3,\k)$ to
$\Omega^1_c(\mathbb{R}^3,\g)$ in the natural way,
\begin{equation}\label{prodc}
<\omega_1,\omega_2>=-\int_{\mathbb{R}^{3}}(\omega_1\wedge,*\overline{\omega}_2)=
-\int_{\mathbb{R}^{3}}*(\omega_1\wedge,*\overline{\omega}_2)d^{3}x,~
\omega_{1,2} \in \Omega^*(\mathbb{R}^{3},\g),
\end{equation}
where $\overline{\omega}_2$ is the complex conjugate of $\omega_2$ with respect to the
complex structure induced by the decomposition $\g=\k\stackrel{\cdot}{+}i\k$.
Note that $\k\subset \g$ is a real subspace with respect to this complex
structure.
We also denote by $\mathcal{H}^0$ the completion of the space
$\Omega^0_c(\mathbb{R}^3,\g)$ with respect to scalar product (\ref{prodc}).

For any $A\in \mathcal{D}$ the operator $\mathrm{div}_A:\Omega^1(\mathbb{R}^{3},\k)\rightarrow
\Omega^{0}(\mathbb{R}^{3},\k)$ gives rise to a linear operator
${\rm div}_A:\Omega^1_c(\mathbb{R}^3,\g)\rightarrow
\Omega^1_c(\mathbb{R}^3,\g)$ and the closure of this operator is a well-defined closed operator
${\rm div}_A:\mathcal{H}^1\rightarrow \mathcal{H}^0$. Denote
by $P_A:\mathcal{H}^1 \rightarrow {\rm Ker}~{\rm div}_A$
the orthogonal projector onto the kernel of this operator. Note that both
$\mathrm{div}_A$ and $P_A$ respect the natural real structure on
$\mathcal{H}^1$, i.e. they are real operators. Let $({\rm Ker}~{\rm div}_A)_c=
{\rm Ker}~{\rm div}_A\cap \Omega^1_c(\mathbb{R}^3,\g)$.

For any $A\in \mathcal{D}$ we define the operator $X_A:({\rm Ker}~{\rm div}_A)_c\rightarrow
{\rm Ker}~{\rm div}_A$ by
\begin{equation}\label{defxa}
X_A=P_A{\rm curl}_A:({\rm Ker}~{\rm div}_A)_c \rightarrow
{\rm Ker}~{\rm div}_A,
\end{equation}
where
\begin{equation}\label{defrot}
{\rm curl}_A=*d_A:({\rm Ker}~{\rm div}_A)_c \rightarrow \Omega^1_c(\mathbb{R}^3,\g).
\end{equation}
\begin{lemma}\label{lemxa}
The operator $X_A$ is real and symmetric in the sense that
$$
<X_Ap,q>=<p,X_Aq>,~\forall p,q\in ({\rm Ker}~{\rm div}_A)_c,
$$
and
$$
\overline{X_Ap}=X_A\overline{p}.
$$
\end{lemma}

\pr

The operator ${\rm curl}_A=*d_A:\Omega^1_c(\mathbb{R}^{3},\g) \rightarrow
\Omega^1_c(\mathbb{R}^{3},\g)$ is a noncommutative analog of the usual $\rm curl$.
As in the commutative case the main property of this operator is that it is
symmetric with respect to scalar product (\ref{prodc}),
$$
<\omega_1,{\rm curl}_A\omega_2>=<{\rm
curl}_A\omega_1,\omega_2>,~\omega_{1,2}\in \Omega^1_c(\mathbb{R}^{3},\g).
$$
This can be checked directly using the Stokes formula and the fact that the
Killing form of $\g$ is invariant with respect to the adjoint action of $\k$ on $\g$.

Since the operator ${\rm curl}_A$ is symmetric $X_A$ is also a symmetric
operator. This completes the proof of the lemma.

\qed

The operators ${\rm curl}_A$, $X_A$ and ${\rm div}_A$ are gauge equivariant in the sense that
$$
{\rm curl}_{g\circ A}=g {\rm curl}_A g^{-1},
$$
$$
X_{g\circ A}=gX_Ag^{-1},
$$
$$
{\rm div}_{g\circ A}=g{\rm div}_Ag^{-1}.
$$
In these formulas it is assumed that the action of the gauge group $\mathcal K$ on the space
$\Omega^*_c(\mathbb{R}^{3},\g)$ is induced by the adjoint representation,
$g\omega=\mathrm{Ad}g(\omega)$, $g\in {\mathcal K}$, $\omega \in \Omega^*_c(\mathbb{R}^{3},\g)$.
Since the Killing form is invariant under the adjoint action the gauge group $\mathcal K$ also acts
 on the spaces $\mathcal{H}^0$ and $\mathcal{H}^1$ by unitary transformations.

We shall prove that the operators
$X_A:({\rm Ker}~{\rm div}_A)_c \rightarrow {\rm Ker}~{\rm div}_A$
have
selfadjoint extensions $X_A:{\rm Ker}~{\rm div}_A\rightarrow
{\rm Ker}~{\rm div}_A$ for $A$ close enough to $0$ with respect to a
certain topology on $\mathcal{D}$, and those selfadjoint extensions are
unitary equivalent to each other.
We start realizing this program by recalling scattering theory for selfadjoint operators
in the form suitable for our purposes.


\section{Wave operators and unitary equivalence}\label{scatt}

Let $\H$ be a complex Hilbert space. In this section we recall, following
\cite{Kato1}, some results on unitary equivalence for selfadjoint operators on $\H$.
Let $\mathcal{C}_0(\H)$ be the set of closed densely defined linear operators $T$ on
$\H$ with domain $\mathfrak D (T)$ and range $\Re (T)$, $\mathcal{B}(\H)$ the set of bounded
operators on $\H$.

Let $T_0$ be a selfadjoint operator in $\H$ and $U\in \mathcal{C}_0(\H)$.
Denote by $R_0(\zeta)=(T_0-\zeta)^{-1}$ the resolvent of $T_0$. As a
function of $\zeta$ the resolvent $R_0(\zeta)$ is holomorphic in the open
lower and upper half-plane.

Suppose that $\mathfrak D(U)\supset \mathfrak D(T_0)$. Then $UR_0(\zeta)\in
\mathcal{B}(\H)$ for ${\rm Im}~\zeta\neq 0$, and $UR_0(\zeta)$ is holomorphic
in the open lower and upper half-plane. The operator $U$ is said to be
$T_0$-smooth if for each $u\in \H$ there is a constant $M_u$ independent of
$\varepsilon$ such that
$$
\int_{-\infty}^{\infty}\| UR_0(\lambda\pm i\varepsilon)u\|^2d\lambda\leq
M_u^2, \forall~\varepsilon>0.
$$
The following proposition gives a strong version of the well-known construction
for the wave operators.
\begin{proposition}\label{uniteq}{\bf (\cite{RS4}, Theorems XIII.24, XIII.26)}
Let $T_0,T$ be two selfadjoint operators in Hilbert space
$\H$ such that $T=T_0+\sum_{i=1}^{n}V^*_iU_i$, where
$U_i,V_i\in \mathcal{C}_0(\H)$ and
$U_i,V_i$ are $T_0$-smooth. Denote by $R_0(\zeta)$ the resolvent of $T_0$ and assume
that
\begin{equation}\label{arb}
A_{ij}=\sup_{\zeta \not\in \mathbb{R}}\|U_iR_0(\zeta)V^*_j\|<\infty,~~i,j=1,\ldots,
n.
\end{equation}
Assume, furthermore, that the norm of the matrix $A=\{A_{ij}\}_{i,j=1,\ldots, n}$,
regarded as an operator in $\mathbb{C}^n$ with the natural Hilbert space norm,
is strictly less than 1,
\begin{equation}\label{uvbound}
\|A\|<1.
\end{equation}

Then the wave operators
$$
W_{\pm}=s-\lim_{t\rightarrow \pm \infty}e^{iTt}e^{-iT_0t}
$$
for the pair $T_0,T$ exist and are unitary
operators, i.e. $W_{\pm}^{-1}\in \mathcal{B}(\H)$. In particular,
the operators $T_0$ and $T$ are
unitary equivalent,
$$
T=W_{\pm}T_0W_{\pm}^{-1}.
$$
\end{proposition}

The condition of $T_0$-smoothness is usually difficult to verify. The
following proposition gives a simple sufficient criterion of $T_0$-smoothness.
\begin{proposition}{\bf (\cite{RS4}, Corollary of Theorem XIII.25)}
Let $T_0$ be a selfadjoint operator in $\H$, $U\in \mathcal{C}_0(\H)$. Then
$U$ is $T_0$-smooth if
\begin{equation}\label{tsmooth}
\sup\|UR_0(\zeta)U^*\|<\infty,
\end{equation}
where the supremum is taken over all $\zeta\in\mathbb{C}$ with ${\rm Im }~
\zeta\neq 0$.
\end{proposition}


\section{Scattering theory for the operators $X_A$: extension to elliptic operators}\label{scattxa}

The main obstruction to direct application of Proposition \ref{uniteq} to
any selfadjoint extensions of
the operators $X_A$ is that for different $A$ these operators act in different
spaces, $X_A:({\rm Ker}~{\rm div}_A)_c\rightarrow
{\rm Ker}~{\rm div}_A$. However, using formula
(\ref{defxa}) one can naturally extend these operators to symmetric operators acting
in $\mathcal H^2$ with the domains $({\rm Ker}~{\rm div}_A)_c\oplus
({\rm Ker}~{\rm div}_A)^\perp_c$, where
$({\rm Ker}~{\rm div}_A)^\perp_c=({\rm Ker}~{\rm div}_A)^\perp\cap \Omega^1_c(\mathbb{R}^{3},\g)$.
If we denote these extensions also by $X_A$ then
\begin{equation}\label{defxa1}
X_A=P_A{\rm curl}_AP_A,
\end{equation}
i.e., we extend $X_A: ({\rm Ker}~{\rm div}_A)_c \rightarrow
({\rm Ker}~{\rm div}_A)$ by zero to the orthogonal complement
$({\rm Ker}~{\rm div}_A)_c$ of its original domain.
Therefore the operators defined by formula (\ref{defxa1}) have big kernels.

As Remark 1.9 in \cite{Kato1} shows, in order to establish
unitary equivalence for operators with nonempty point spectrum we
have to significantly restrict the class of perturbations (actually if $u$ is an eigenvector
of $T_0$, in the notation of Proposition \ref{uniteq}, then for $U$ to be $T_0$-smooth
we must have $Uu=0$). This is not satisfactory for our purposes.

In order to overcome this difficulty we shall find an elliptic extension $S_A$
for the operator $X_A$ (compare with \cite{BS1,BS2}). $S_A$ is a symmetric operator acting in
the space
$\mathcal H^1\stackrel{\cdot}{+}
\mathcal H^0$ with the domain $\mathfrak D (S_A)=(({\rm Ker}~{\rm div}_A)_c\oplus
({\rm Ker}~{\rm div}_A)^\perp_c)\stackrel{\cdot}{+} \Omega^0_c(\mathbb{R}^{3},\g)$.
If we represent $S_A$
in the block form then
\begin{equation}\label{sa}
S_A=\left(
\begin{array}{cc}
P_A{\rm curl}_AP_A & -d_A \\
  {\rm div}_A & 0
\end{array}
\right).
\end{equation}

In Section \ref{scattxa1} we shall define natural selfadjoint extensions for $S_A$
and obtain wave operators for selfadjoint extensions of operators
$X_A$ using those for selfadjoint extensions of $S_A$.


\section{The properties of the unperturbed operator $S_0$}\label{S0}

In order to apply Proposition \ref{uniteq} to the operators $S_A$ we have to
study first the unperturbed operator $S_0$. In particular, in view of
condition (\ref{arb}) we have to study the properties of the resolvent of
the operator $S_0$.

First we note that the operator $S_0$,
\begin{equation}\label{s0}
S_0=\left(
\begin{array}{cc}
{\rm curl} & -d \\
  {\rm div} & 0
\end{array}
\right),
\end{equation}
acting in the space $\mathcal H^1\stackrel{\cdot}{+}\mathcal H^0$ with the
natural domain
$$
\mathfrak D (S_0)=\{(\omega,u)\in\mathcal H^1\stackrel{\cdot}{+}\mathcal H^0:{\rm curl}~\omega,
du\in\mathcal{H}^1,{\rm div}~\omega\in \mathcal{H}^0\}
$$
is selfadjoint (see \cite{BS1,BS2}). The formula for the resolvent of this
operator can easily be obtained with the help of an explicit expression for $S_0^2$,
\begin{equation}\label{s02}
S_0^2=\left(
\begin{array}{cc}
\widehat{\triangle} & 0 \\
  0 & \triangle
\end{array}
\right),
\end{equation}
where $\triangle$ is the usual Laplace operator and $\widehat{\triangle}$
stands for the Laplace operator acting on the components of elements from
$\mathcal{H}^1$. If we denote the resolvent of $S_0$ by $R_0(\lambda)$ then
(see \cite{Fi} for a similar calculation for the operator $\rm curl$)
$$
R_0(\lambda)=(S_0-\lambda I)^{-1}=(S_0+\lambda I)(S_0^2-\lambda^2 I)^{-1},
$$
or in components
\begin{equation}\label{rs0}
\begin{array}{lr}
R_0(\lambda)=\left(
\begin{array}{cc}
{\rm curl}+\lambda I & -d \\
  {\rm div} & \lambda I
\end{array}
\right)
\left(
\begin{array}{cc}
(\widehat{\triangle}-\lambda^2 I)^{-1} & 0 \\
  0 & (\triangle-\lambda^2 I)^{-1}
\end{array}
\right)= & \\
 & \\
& \\ \hspace{3cm} \left(
\begin{array}{cc}
({\rm curl}+\lambda I)(\widehat{\triangle}-\lambda^2 I)^{-1} & -d(\triangle-\lambda^2 I)^{-1} \\
  {\rm div}(\widehat{\triangle}-\lambda^2 I)^{-1} & \lambda(\triangle-\lambda^2 I)^{-1}
\end{array}
\right).
\end{array}
\end{equation}

Verification of condition (\ref{arb}) will be based on the fact that the
resolvent $R_0(\lambda)$ acts as a bounded operator in certain weighted $L^2$-spaces.

We recall that for any real $s$ the weighted space $L^{2,s}(\mathbb{R}^{n})$
is defined by
$$
L^{2,s}(\mathbb{R}^{n})=\{u(x):~(1+|x|^2)^{\frac s2}u(x)\in
L^{2}(\mathbb{R}^{n})\},
$$
where $|x|$ is the usual norm of $x$ in $\mathbb{R}^{n}$. The space
$L^{2,s}(\mathbb{R}^{n})$ is a Hilbert space, the corresponding norm $\|\cdot \|_s$ is
equal to
$$
\|u\|_s=\|(1+|x|^2)^{\frac s2}u(x)\|_{L^{2}(\mathbb{R}^{n})}.
$$
Clearly, the spaces $L^{2,s}(\mathbb{R}^{n})$ and $L^{2,-s}(\mathbb{R}^{n})$
are dual to each other.

We shall denote by $\mathcal{H}^{0,s}$ and $\mathcal{H}^{1,s}$ the weighted
versions of the spaces $\mathcal{H}^0$ and $\mathcal{H}^1$, with the
obvious modifications of the scalar products.

For any bounded operator $T:L^{2,s_1}(\mathbb{R}^{n})\rightarrow L^{2,s_2}(\mathbb{R}^{n})$
we denote by $\|T\|_{s_1,s_2}$ its norm. The required estimates
for the resolvent $R_0$ will be based on the following well-known results by Agmon
and Jensen-Kato on the resolvent of the Laplace operator.
\begin{proposition}{\bf (\cite{Ag}, Appendix A, Remark 2; \cite{JK}, Lemma
2.1)}\label{reslapl}
Let $(\triangle-\zeta I)^{-1}$ be the resolvent of the Laplace operator acting in
$L^2(\mathbb{R}^{3})$. Then for any $\lambda \in \mathbb{C}$ and $s>1$ the operators
\begin{equation}\label{resblocks}
\lambda(\triangle-\lambda^2I)^{-1},~~\partial_i(\triangle-\lambda^2
I)^{-1},~~i=1,2,3,
\end{equation}
act as bounded operators from $L^{2,s}(\mathbb{R}^{n})$ to
$L^{2,-s}(\mathbb{R}^{n})$. Moreover,
$$
\begin{array}{l}
\|\lambda(\triangle-\lambda^2I)^{-1}\|_{s,-s}\leq C, \\
\\
\|\partial_i(\triangle-\lambda^2I)^{-1}\|_{s,-s}\leq C,~~i=1,2,3,
\end{array}
$$
where $C$ is a constant independent of $\lambda$.

Fix $\varepsilon >0$. Then for any $\lambda \in \mathbb{C}$ with $|\lambda|\geq \varepsilon$
and $s>\frac{1}{2}$ the operators
$$
\frac{1}{\lambda}\partial_i\partial_j(\triangle-\lambda^2
I)^{-1},~~i,j=1,2,3
$$
act as bounded operators from $L^{2,s}(\mathbb{R}^{n})$ to
$L^{2,-s}(\mathbb{R}^{n})$. Moreover,
$$
\|\frac{1}{\lambda}\partial_i\partial_j(\triangle-\lambda^2I)^{-1}\|_{s,-s}\leq C,~~i=1,2,3,
$$
where $C$ is a constant independent of $\lambda$.
\end{proposition}

Clearly, the blocks of the resolvent $R_0(\lambda)$ consist of the terms of
 form (\ref{resblocks}) (see formula (\ref{rs0})). Therefore Proposition
\ref{reslapl} provides an estimate uniform in $\lambda$ of the norm of the resolvent
$R_0(\lambda)$ acting as a bounded operator from $\mathcal{H}^{1,s}
\stackrel{\cdot}{+}\mathcal H^{0,s}$ to
$\mathcal{H}^{1,-s}\stackrel{\cdot}{+}\mathcal H^{0,-s}$. We formulate this result as a
corollary to Proposition \ref{reslapl}.
\begin{corollary}\label{robound}
For any $\lambda \in \mathbb{C}$ and $s>1$ the resolvent
$R_0(\lambda)$ acts as a bounded operator from $\mathcal{H}^{1,s}
\stackrel{\cdot}{+}\mathcal H^{0,s}$ to
$\mathcal{H}^{1,-s}\stackrel{\cdot}{+}\mathcal H^{0,-s}$. Moreover,
$$
\|R_0(\lambda)\|_{s,-s}\leq C,
$$
where $C$ is a constant independent of $\lambda$.
\end{corollary}


\section{Scattering theory for the elliptic operators $S_A$}\label{scattsa}

In this section we find conditions under which the operators $S_0$ and $S_A$ are unitary
equivalent. In order to derive these conditions we apply Proposition
\ref{uniteq} to the pair $S_0,S_A$. First we formulate our main result.
\begin{theorem}\label{sauniteq}
Let $A$ be a $\k$-valued connection 1-form on $\mathbb{R}^{3}$,
$A\in \Omega^1(\mathbb{R}^{3},\k)$, $F$ the curvature of $A$.
There exists a positive constant $C$ such that if for some $s>3$
and $g\in \mathcal{K}$
\begin{eqnarray}
\|g\circ F(1+|x|^2)^{s}\|_{(\infty)}< C, \label{fbound} \\
\|g\circ A(1+|x|^2)^{\frac{s}{2}}\|_{(\infty)}< C, \label{abound}
\end{eqnarray}
where for any $\omega \in \Omega^i(\mathbb{R}^{3},\g)$
$$
\|\omega \|_{(\infty )}={\rm ess}-\sup_{x\in
\mathbb{R}^3}*(\omega(x)\wedge,*\overline{\omega(x)})^{\frac{1}{2}},
$$
then

(i) The corresponding operator $S_A$ is selfadjoint with the domain $g^{-1}\mathfrak D (S_0)$.

(ii) The wave operators
$$
W_{\pm}(S_0,S_A)=s-\lim_{t\rightarrow \pm \infty}e^{iS_At}e^{-iS_0t}
$$
for the pair $S_0,S_A$ exist and are unitary. In particular,
the selfadjoint operators $S_0$ and $S_A$ are unitary equivalent.
\end{theorem}
\begin{remark}
Since the operators $S_A$ are gauge equivariant and the gauge group $\mathcal K$ acts on the space
$\mathcal{H}^1$ by unitary transformations we shall assume,
without loss of generality, that conditions (\ref{fbound}),
(\ref{abound}) are imposed on $F$ and $A$, or, in other words,
that $A$ is in the gauge in which the conditions (\ref{fbound}),
(\ref{abound}) are satisfied.
\end{remark}

We start the proof of this theorem with the study
of the orthogonal projection operator $P_A$ and the ``magnetic'' Laplace
operator $\triangle_A=-{\rm div}_Ad_A$.
\begin{lemma}\label{deltaa}
Let $\mathcal{H}^1_{loc}(\mathbb{R}^3,\k)$ be the space of locally square integrable,
with respect to the
scalar product (\ref{prod}), $\k$-valued 1-forms on $\mathbb{R}^3$. Then for any
$A\in \mathcal{H}^1_{loc}(\mathbb{R}^3,\k)$ we have:

(i) The operators $-d_A:\Omega^0_c(\mathbb{R}^{3},\g)
\rightarrow \mathcal{H}^1$ and ${\rm div}_A:\Omega^1_c(\mathbb{R}^{3},\g)\rightarrow
\mathcal{H}^0$ are closable and their closures are operators
$-d_A:\mathcal{H}^0\rightarrow
\mathcal{H}^1$, ${\rm div}_A:\mathcal{H}^1\rightarrow \mathcal{H}^0$ which
are adjoint to each other.

(ii) The operator $\triangle_A: \mathcal{H}^0\rightarrow
\mathcal{H}^0$ defined by the differential expression $\triangle_A=-{\rm
div}_Ad_A$ is selfadjoint on the natural domain $\mathfrak D(\triangle_A)=
\{u\in \mathfrak D (d_A)|~d_Au\in \mathfrak D ({\rm div}_A)\}$.

(iii) The operator $\triangle_A$ has trivial kernel, and the inverse operator $\triangle_A^{-1}$
is well-defined.

(iv) The operator $P_A$ is the closure of the operator $I+d_A\triangle_A^{-1}{\rm
div}_A$ defined on $\Re (d_A)\bigoplus {\rm Ker}~{\rm div}_A$.
\end{lemma}

\pr
For parts (i) and (ii) see \cite{RS2}, \S X.3, Example 4.

In order to prove part (iii) we observe that if $u\in {\rm Ker}~\triangle_A$ then $d_Au=0$
since $0=(\triangle_Au,u)=-({\rm div}_Ad_Au,u)=(d_Au,d_Au)$. But then by the
invariance of the Killing form we have
$$
d(u,u)=(d_Au,u)+(u,d_Au)=0.
$$
Therefore $(u,u)={\rm const}$, and $u$ does not belong to $\mathcal{H}^0$.

For part (iv) we first note that the orthogonal direct sum
$\Re (d_A)\bigoplus {\rm Ker}~{\rm div}_A$ is dense
in $\mathcal{H}^1$, $\mathcal{H}^1=\overline{\Re (d_A)}\bigoplus
{\rm Ker}~{\rm div}_A$

Next, if $\omega \in \Re (d_A),~\omega =d_Au$ then
$(I+d_A\triangle_A^{-1}{\rm div}_A)\omega=\omega+d_A\triangle_A^{-1}{\rm
div}_Ad_Au=0$ since $\triangle_A^{-1}{\rm div}_Ad_Au=-u$, and if $\omega \in
{\rm Ker}~{\rm div}_A$ then $(I+d_A\triangle_A^{-1}{\rm
div}_A)\omega=\omega$. This completes the proof.

\qed

Now we write down the perturbation $W=S_A-S_0$ in a
convenient form. From formula (\ref{sa}) we formally have
\begin{equation}\label{perturb}
W=\left(
\begin{array}{cc}
P_0*\mathrm{ad}A P_0-P_0{\rm curl}_A\Delta P -\Delta P{\rm curl}_AP_0 +
\Delta P{\rm curl}_A\Delta P& -\mathrm{ad}A \\
 *\mathrm{ad}A* & 0
\end{array}
\right) ,
\end{equation}
where $\mathrm{ad}A(\omega)=[A\wedge \omega]$, $\omega \in
\Omega^i(\mathbb{R}^3,\g)$, and $\Delta P=P_0-P_A$.

Using the expression for the operator $P_A$ obtained in part (iv) of Lemma
\ref{deltaa}, the Hilbert-type identity,
$$
\triangle_A^{-1}=\triangle^{-1}-\triangle^{-1}(\triangle_A-\triangle)\triangle_A^{-1},
$$
and the formula
$$
\triangle_A-\triangle=-{\rm
div}~\mathrm{ad}A-*\mathrm{ad}A*d-*\mathrm{ad}A*\mathrm{ad}A
$$
one can derive the following formal expression for $\Delta P$:
\begin{eqnarray}\label{deltap}
\Delta P=d\triangle^{-1}*\mathrm{ad}A*+\mathrm{ad}A\triangle^{-1}\mathrm{div}_A+ \hspace{5cm} \\
\hspace{3cm} +d_A\triangle^{-1}\mathrm{div}~\mathrm{ad}A\triangle_A^{-1}\mathrm{div}_A+
d_A\triangle^{-1}*\mathrm{ad}A*d_A\triangle_A^{-1}\mathrm{div}_A. \nonumber
\end{eqnarray}

Substituting (\ref{deltap}) into formula (\ref{perturb}) we get useful
expressions for the terms of the perturbation $W$,
\begin{eqnarray}
\Delta P{\rm
curl}_AP_0=\left(
d\triangle^{-1}*\mathrm{ad}F-d\triangle^{-1}\mathrm{div}*\mathrm{ad}A+
\mathrm{ad}A\triangle^{-1}*\mathrm{ad}F+ \right. \hspace{1cm} \nonumber \\
\hspace{1.5cm} +d\triangle^{-1}\mathrm{div}~\mathrm{ad}A\triangle_A^{-1}*\mathrm{ad}F+
d\triangle^{-1}*\mathrm{ad}A*d_A\triangle_A^{-1}*\mathrm{ad}F+ \label{dprotapo}\\
\hspace{1.5cm} \left.
+\mathrm{ad}A\triangle^{-1}\mathrm{div}~\mathrm{ad}A\triangle_A^{-1}*\mathrm{ad}F+
\mathrm{ad}A\triangle^{-1}*\mathrm{ad}A*d_A\triangle_A^{-1}*\mathrm{ad}F
\right) P_0, \nonumber
\end{eqnarray}
\begin{eqnarray}
P_0{\rm curl}_A\Delta P=P_0\left(*\mathrm{ad}F\triangle^{-1}\mathrm{div}-
*\mathrm{ad}Ad\triangle^{-1}\mathrm{div}
+*\mathrm{ad}F\triangle^{-1}*\mathrm{ad}A*+ \right. \nonumber \\
\hspace{1cm} +*\mathrm{ad}F\triangle^{-1}_A*\mathrm{ad}A*d\triangle^{-1}\mathrm{div}
+*\mathrm{ad}F\triangle^{-1}_A\mathrm{div}_A\mathrm{ad}A\triangle^{-1}\mathrm{div}+ \label{porotadp}\\
\hspace{1cm} \left. +*\mathrm{ad}F\triangle^{-1}_A*\mathrm{ad}A*d\triangle^{-1}*\mathrm{ad}A*
+*\mathrm{ad}F\triangle^{-1}_A\mathrm{div}_A\mathrm{ad}A\triangle^{-1}*\mathrm{ad}A*
\right), \nonumber
\end{eqnarray}
\begin{eqnarray}
\Delta P{\rm curl}_A\Delta P=\left(d\triangle^{-1}*\mathrm{ad}A*+
\mathrm{ad}A\triangle^{-1}\mathrm{div}+\mathrm{ad}A\triangle^{-1}*\mathrm{ad}A*+ \right. \hspace{1cm} \nonumber \\
+d\triangle^{-1}\mathrm{div}~\mathrm{ad}A\triangle^{-1}_A\mathrm{div}_A
+\mathrm{ad}A\triangle^{-1}\mathrm{div}~\mathrm{ad}A\triangle^{-1}_A\mathrm{div}_A+ \nonumber \\
\left. +d\triangle^{-1}*\mathrm{ad}A*d_A\triangle^{-1}_A\mathrm{div}_A
+\mathrm{ad}A\triangle^{-1}*\mathrm{ad}A*d_A\triangle^{-1}_A\mathrm{div}_A
\right)\times \label{dprotadp}\\
\times \left( *\mathrm{ad}F\triangle^{-1}\mathrm{div}-
*\mathrm{ad}Ad\triangle^{-1}\mathrm{div}
+*\mathrm{ad}F\triangle^{-1}*\mathrm{ad}A*+ \right. \hspace{3cm} \nonumber \\
\hspace{1cm} +*\mathrm{ad}F\triangle^{-1}_A*\mathrm{ad}A*d\triangle^{-1}\mathrm{div}
+*\mathrm{ad}F\triangle^{-1}_A\mathrm{div}_A\mathrm{ad}A\triangle^{-1}\mathrm{div}+ \nonumber \\
\hspace{1cm} \left. +*\mathrm{ad}F\triangle^{-1}_A*\mathrm{ad}A*d\triangle^{-1}*\mathrm{ad}A*
+*\mathrm{ad}F\triangle^{-1}_A\mathrm{div}_A\mathrm{ad}A\triangle^{-1}*\mathrm{ad}A*
\right). \nonumber
\end{eqnarray}

In order to analyze expressions (\ref{dprotapo})--(\ref{dprotadp}) we shall need the properties
of $\triangle_A^{-1}$ as an operator acting in
weighted $L^2$-spaces and $L^p$-spaces. For any $p\geq 1$ denote by $\mathcal{H}^i_p$
the closure of the space
$\Omega^i_c(\mathbb{R}^{3},\g)$ with respect to the norm
$$
\|\omega \|_{(p)}=\left( -\int_{\mathbb{R}^{3}}*(\omega\wedge,*\overline{\omega})^{\frac{p}{2}}d^{3}x
\right)^{\frac{1}{p}},~p<\infty,
$$
$$
\|\omega \|_{(\infty )}={\rm ess-}\sup_{x\in
\mathbb{R}^3}*(\omega(x)\wedge,*\overline{\omega(x)})^{\frac{1}{2}}.
$$

Note that from the H\"{o}lder inequality it follows that for any $p\geq 2$ we have a
natural embedding, $\mathcal{H}^i_p\subset \mathcal{H}^{i,-s}$, where $s>3(\frac{1}{2}-
\frac{1}{p})$, and for any $1\leq p\leq 2$ we have another embedding, $\mathcal{H}^{i,s}\subset
\mathcal{H}^i_p$, where $s>3(\frac{2}{p}-1)$, i.e.,
\begin{equation}\label{psimbed}
\begin{array}{l}
\mathcal{H}^i_p\subset \mathcal{H}^{i,-s},~p\geq 2,~s>3(\frac{1}{2}-
\frac{1}{p}), \\
 \\
\mathcal{H}^{i,s}\subset \mathcal{H}^i_p,~1\leq p\leq 2,~s>3(\frac{2}{p}-1).
\end{array}
\end{equation}

Obviously, there are also natural embeddings
\begin{equation}\label{ssimbed}
\mathcal{H}^{i,s_1}\subset \mathcal{H}^{i,s_2}\subset \mathcal{H}^{i}\subset
\mathcal{H}^{i,-s_2}\subset \mathcal{H}^{i,-s_1},~s_1>s_2>0.
\end{equation}

\begin{lemma}\label{daact}
(i) For any
$A\in \mathcal{H}^1_{loc}(\mathbb{R}^3,\k)$, $u\in \mathcal{H}^0_{\frac{6}{5}}$ and
$\omega \in \mathcal{H}^1_{2}$
we have
\begin{eqnarray}
\|\triangle_A^{-1}u\|_{(6)}\leq K^2\|u\|_{(\frac{6}{5})}, \label{dabound} \\
\|d_A\triangle_A^{-1}u\|_{(2)}\leq K\|u\|_{(\frac{6}{5})}, \label{dadabound} \\
\|\triangle_A^{-1}{\rm div}_A\omega \|_{(6)}\leq K\|\omega \|_{(2)}, \label{dadivbound}
\end{eqnarray}
where $K$ is a constant independent of $A$.

(ii) The following operators are bounded and have norms independent of
$A$:
\begin{equation}\label{dabounds}
\begin{array}{l}
\triangle_A^{-1}:\mathcal{H}^{0,s_1}\rightarrow \mathcal{H}^{0,-s_2},~s_1>2,~s_2>1,\\
\\
\triangle_A^{-1}:\mathcal{H}^{0,s_1}\rightarrow \mathcal{H}^{0,-s_2},~s_1>1,s_2>2\\
\\
\triangle_A^{-1}{\rm div}_A:\mathcal{H}^{1}\rightarrow \mathcal{H}^{0,-s},~s>1\\
\\
d_A\triangle_A^{-1}:\mathcal{H}^{0,s}\rightarrow \mathcal{H}^{1},~s>1.
\end{array}
\end{equation}
\end{lemma}

\pr
Inequality (\ref{dabound}) follows from the well-known Kato inequality,
$\|d|u|\|_{(2)}\leq \|d_Au\|_{(2)}$ (see \cite{JT}, Appendix), and the isoperimetric-type
inequality, $\|f\|_{(6)}\leq K\|df\|_{(2)}$, $f\in C_0^\infty
(\mathbb{R}^{3})$ (see \cite{LU}, Theorem 2.1). Indeed, these two
inequalities and the H\"{o}lder inequality imply that
\begin{equation}\label{dabound1}
\|w\|_{(6)}^2\leq K^2\|d|w|\|_{(2)}^2\leq K^2(\triangle_Aw,w)\leq K^2\|w\|_{(6)}
\|\triangle_A w\|_{(\frac{6}{5})}.
\end{equation}
This inequality holds for $w\in \mathcal{C}$, where $\mathcal{C}=\{\varphi v|
\varphi \in C_0^\infty
(\mathbb{R}^{3}),~v\in (\triangle_A+1)^{-1}(\mathcal{H}^0_2\bigcap
\mathcal{H}^0_\infty)\}$ is an operator core of $\triangle_A$ (see
\cite{LS}, Lemma 5).

Now dividing both sides of (\ref{dabound1}) by $\|w\|_{(6)}$ and denoting
$u=\triangle_Aw$ we get (\ref{dabound}).

Similarly, for any $u=\triangle_Aw$, $w\in \mathcal{C}$ inequality (\ref{dabound}) and the
H\"{o}lder inequality yield
$$
\|d_A\triangle_A^{-1}u\|_{(2)}^2=(\triangle_A\triangle_A^{-1}u,\triangle_A^{-1}u)=
(u,\triangle_A^{-1}u)\leq
\|\triangle_A^{-1}u\|_{(6)}\|u\|_{(\frac{6}{5})}\leq K^2
\|u\|_{(\frac{6}{5})}^2,
$$
which is equivalent to (\ref{dadabound}).

Now we deduce (\ref{dadivbound}) from (\ref{dadabound}) by duality.

Finally we infer part (ii) of the lemma using part (i), embeddings (\ref{psimbed}) and duality
$(\mathcal{H}^{i,s})^*\simeq \mathcal{H}^{i,-s}$.

\qed

{\em Proof of part (i) of Theorem \ref{sauniteq}.}
First we note that conditions (\ref{fbound}) and (\ref{abound}) ensure that
for any real $s'$
\begin{equation}\label{faoper}
\begin{array}{l}
\mathrm{ad}F: \mathcal{H}^{i,s'}\rightarrow \mathcal{H}^{i+2,s'+2s}, \\
\\
\mathrm{ad}A: \mathcal{H}^{i,s'}\rightarrow \mathcal{H}^{i+1,s'+s}\\
\end{array}
\end{equation}
are bounded operators with norms not exceeding $C$. Combining this fact with
part (ii) of Lemma \ref{daact} and recalling embeddings (\ref{ssimbed}) we
infer that the operators $*\mathrm{ad}A:\mathcal{H}^1\rightarrow \mathcal{H}^1$,
 $\mathrm{ad}A:\mathcal{H}^0\rightarrow \mathcal{H}^1$,
$*\mathrm{ad}A*:\mathcal{H}^1\rightarrow \mathcal{H}^0$ are bounded and the operators defined
by formulas (\ref{dprotapo})--(\ref{dprotadp}) are bounded in the space
$\mathcal{H}^1$. Therefore the perturbation $W$ is a bounded operator, and
part (i) of  Theorem \ref{sauniteq} follows from the Kato-Rellich theorem (see
\cite{RS2}, Theorem X.12).

Now we discuss part (ii) of Theorem \ref{sauniteq}. In order to apply
Theorem \ref{uniteq} to the pair of operators $S_0, S_A$ we have to
factorize every term of the perturbation $W$ and then verify conditions
(\ref{uvbound}), (\ref{tsmooth}). Note that the perturbation $W$ contains not
only multiplication operators but also integral operators (see formulas
(\ref{dprotapo})--(\ref{dprotadp})) and to check conditions (\ref{uvbound}), (\ref{tsmooth})
we shall need not only estimates for the norm of the resolvent
$R_0(\lambda)$ of the unperturbed operator $S_0$ obtained in Proposition \ref{reslapl}
and Corollary \ref{robound} but also estimates of norms of compositions of
$R_0(\lambda)$ and of some integral operators. More precisely, in view of
(\ref{dprotapo})--(\ref{dprotadp}), we have to obtain the following
estimates.
\begin{lemma}\label{rointbounds}
For every $\lambda \in \mathbb{C}$ the operators
\begin{eqnarray}
-\triangle^{-1}\mathrm{div}(\mathrm{curl}+\lambda
I)(\widehat{\triangle}-\lambda^{2}I)^{-1}d\triangle^{-1}= \label{ddrdd} \hspace{3cm} \\
\hspace{3cm} =\lambda\triangle^{-1}
(\triangle-\lambda^{2}I)^{-1}:\mathcal{H}^{0,s_1}\rightarrow
\mathcal{H}^{0,-s_2},~s_1>3, s_2>\frac{3}{2}, \nonumber \\
(\mathrm{curl}+\lambda
I)(\widehat{\triangle}-\lambda^{2}I)^{-1}d\triangle^{-1}= \label{rdd} \hspace{4.5cm} \\
\hspace{3cm} =\lambda d\triangle^{-1}
(\triangle-\lambda^{2}I)^{-1}:\mathcal{H}^{0,s}\rightarrow
\mathcal{H}^{1,-s},~s>1, \nonumber \\
\triangle^{-1}\mathrm{div}(\mathrm{curl}+\lambda
I)(\widehat{\triangle}-\lambda^{2}I)^{-1}= \label{ddr} \hspace{4cm} \\
\hspace{3cm} =\lambda
\mathrm{div}\widehat{\triangle}^{-1}
(\widehat{\triangle}-\lambda^{2}I)^{-1}:\mathcal{H}^{1,s}\rightarrow
\mathcal{H}^{0,-s},~s>1, \nonumber \\
(\mathrm{curl}+\lambda
I)(\widehat{\triangle}-\lambda^{2}I)^{-1}d\triangle^{-1}\mathrm{div}= \label{rp} \hspace{4cm} \\
\hspace{3cm} =\lambda
(\widehat{\triangle}-\lambda^{2}I)^{-1}d\triangle^{-1}\mathrm{div}:\mathcal{H}^{1,s}
\rightarrow \mathcal{H}^{1,-s},~s>1, \nonumber \\
-d\triangle^{-1}\mathrm{div}(\mathrm{curl}+\lambda
I)(\widehat{\triangle}-\lambda^{2}I)^{-1}d\triangle^{-1}\mathrm{div}= \label{prp} \hspace{2cm} \\
\hspace{3cm} =\lambda
(\widehat{\triangle}-\lambda^{2}I)^{-1}d\triangle^{-1}\mathrm{div}:\mathcal{H}^{1,s}
\rightarrow \mathcal{H}^{1,-s},~s>1. \nonumber
\end{eqnarray}
are bounded and have norms not exceeding $D$, where $D$ is a constant independent of
$\lambda$.
\end{lemma}

\pr
First we recall
that the resolvent of the Laplace operator is an integral operator
explicitly given by the following formula (see \cite{RS2})
\begin{equation}\label{rlapform}
((\triangle-\lambda^{2} I)^{-1}u)(x)=\int_{\mathbb{R}^{3}}
\frac{e^{i\lambda{\rm sgn}({\rm Im}~\lambda)|x-y|}}{4\pi |x-y|}u(y)d^{3}y.
\end{equation}

To prove (\ref{ddrdd}) we shall use the Hilbert identity in the form
\begin{equation}\label{dd}
\lambda\triangle^{-1}(\triangle-\lambda^{2}I)^{-1}=\frac{1}{\lambda}
((\triangle-\lambda^{2}I)^{-1}-\triangle^{-1}).
\end{equation}

From Proposition \ref{reslapl} it follows that for any $\varepsilon >0$ and
$|\lambda|>\varepsilon$ the operator in the r.h.s. of (\ref{dd}) is a
bounded operator acting from $\mathcal{H}^{0,s_1}$ to $\mathcal{H}^{0,-s_2}$ for
any $s_{1,2}>1$, with the norm uniformly bounded in $\lambda$. We have to check that
the norm of this operator remains
finite when $\lambda \rightarrow 0$ and $s_1>3,~s_2>\frac{3}{2}$.

Indeed, using formulas (\ref{rlapform}) and (\ref{dd}) we obtain that for any
$u\in \Omega_c^0(\mathbb{R}^3,\g)$
$$
\lim_{\lambda \rightarrow 0}\lambda\triangle^{-1}(\triangle-\lambda^{2}I)^{-1}u
=\frac{i}{4\pi}\int_{\mathbb{R}^3}u(y)d^3y.
$$
The operator in the r.h.s. of the last formula is bounded from
$\mathcal{H}^0_1$ to $\mathcal{H}^0_{\infty}$, and hence, in view of
embeddings (\ref{psimbed}), from $\mathcal{H}^{0,s_1}$ to $\mathcal{H}^{0,-s_2}$
for $s_1>3,~s_2>\frac{3}{2}$. This proves (\ref{ddrdd}).

(\ref{rdd}) and (\ref{ddr}) can be proved in a similar way with the help of formula (\ref{dd}).
But one should also apply Lemma 2.1 in \cite{JK} and instead of
formula (\ref{rlapform}) one should use the following expressions for the
terms of the kernels of the operators in the r.h.s. of (\ref{rdd}) and
(\ref{ddr}):
\begin{eqnarray}\label{dres}
\partial_i\frac{e^{i\lambda{\rm sgn}({\rm Im}~\lambda)|x-y|}-1}{4\pi
|x-y|\lambda}=\frac{ie^{i\lambda{\rm sgn}({\rm Im}~\lambda)|x-y|}(x_i-y_i)}{4\pi
|x-y|^2} \hspace{1cm} \\
\hspace{6cm} -\frac{(e^{i\lambda{\rm sgn}({\rm Im}~\lambda)|x-y|}-1)(x_i-y_i)}{4\pi
|x-y|^3\lambda}. \nonumber
\end{eqnarray}

The behavior of the norms of the operators (\ref{rp}) and (\ref{prp})
when $\lambda \rightarrow 0$ is controlled with the help of the formula for the kernels of
the operators
$$
\partial_i\partial_j\lambda\triangle^{-1}(\triangle-\lambda^{2}I)^{-1}
$$
that is similar to (\ref{dres}). The fact that these norms are finite when
$\lambda \rightarrow \infty$ can be proved using formula (\ref{dd}) and the
last statement of Proposition \ref{reslapl}.

\qed

{\em Proof of part (ii) of Theorem \ref{sauniteq}.}
We have to factorize every term of the perturbation $W$ and then verify conditions
(\ref{uvbound}), (\ref{tsmooth}). We demonstrate how to obtain the required
estimates in case of the most ``singular'' terms. All the other terms in
(\ref{perturb}) can be
treated in a similar way using expressions (\ref{dprotapo})--(\ref{dprotadp}),
conditions (\ref{fbound}), (\ref{abound}), part (ii)
of Lemma \ref{daact}, Proposition \ref{reslapl},
Corollary \ref{robound} and Lemma \ref{rointbounds}.

First, let us consider formula (\ref{dprotadp}). This term of the
perturbation is already factorized. Let us consider the part
$d\triangle^{-1}*\mathrm{ad}A*$ in the first term of (\ref{dprotadp}). According to
formula (\ref{tsmooth}) one
should check that the composition $\mathrm{ad}A\triangle^{-1}\mathrm{div}(\mathrm{curl}+\lambda
I)(\widehat{\triangle}-\lambda^{2}I)^{-1}d\triangle^{-1}*\mathrm{ad}A*$ is a
bounded operator in $\mathcal{H}^1$ with the norm uniformly bounded in $\lambda$.
This follows from part (\ref{ddrdd}) of
Lemma \ref{rointbounds} and from (\ref{faoper}).

Now let us consider the first term in (\ref{dprotapo}),
$d\triangle^{-1}*\mathrm{ad}FP_0$. If we write $F=\sum_{a,i,j}F_{ij}^aT_adx_i\wedge
dx_j$, where $T_a$ is a basis of $\mathfrak k$, then that term is, in turn, the sum of the
following ones:
$d\triangle^{-1}F_{ij}^a*\mathrm{ad}(T_adx_i\wedge dx_j)P_0$. Since $F_{ij}^a$
are real-valued functions each of these terms can be factorized as follows:
\begin{eqnarray*}
d\triangle^{-1}F_{ij}^a*\mathrm{ad}(T_adx_i\wedge dx_j)P_0= \hspace{7cm}\\
\hspace{3cm} =\left( d\triangle^{-1}\sqrt{|F_{ij}^a|}\right)
\left( \sqrt{|F_{ij}^a|}\mathrm{sgn}(F_{ij}^a)*\mathrm{ad}(T_adx_i\wedge
dx_j)P_0\right) .
\end{eqnarray*}

Now recall conditions (\ref{uvbound}) and (\ref{tsmooth}).
Let us show, for instance, that the composition
$\sqrt{|F_{ij}^a|}\triangle^{-1}\mathrm{div}(\mathrm{curl}+\lambda
I)(\widehat{\triangle}-\lambda^{2}I)^{-1}d\triangle^{-1}\sqrt{|F_{ij}^a|}$
is a bounded operator in $\mathcal{H}^0$ with the norm uniformly bounded in $\lambda$.
This follows from part (\ref{ddrdd}) of
Lemma \ref{rointbounds} and from the fact that by (\ref{fbound}) the
operator of multiplication by $\sqrt{|F_{ij}^a|}$ is bounded from
$\mathcal{H}^{0,s'}$ to $\mathcal{H}^{0,s'+s}$ for any real $s'$.
This completes the proof of Theorem \ref{sauniteq}.

\qed


\section{Scattering theory for the operators $X_A$: reduction from elliptic operators $S_A$}
\label{scattxa1}

In this section we show how to prove that selfadjoint extensions of symmetric operators
$X_A$ defined in Lemma \ref{lemxa} are unitary equivalent. We start with
the definition of these extensions.

Observe that the subspace $\mathrm{Ker}~\mathrm{div}_A$, as well as its orthogonal complement,
is an invariant
subspace in $\mathcal{H}^1\stackrel{\cdot}{+}\mathcal{H}^0$ for the selfadjoint operator
$S_A$ defined in Theorem \ref{sauniteq}. The restriction of $S_A$ to this
subspace gives an operator that is, obviously, a selfadjoint extension of
the symmetric operator $X_A$ introduced in Lemma \ref{lemxa}. We denote this
extension by the same letter, $X_A:\mathrm{Ker}~\mathrm{div}_A\rightarrow
\mathrm{Ker}~\mathrm{div}_A$.
\begin{theorem}\label{xauniteq}
Suppose that the conditions of Theorem \ref{sauniteq} are satisfied.
Then the operators $U_\pm=P_AW_\pm(S_0,S_A)P_0:\mathrm{Ker}~\mathrm{div}\rightarrow
\mathrm{Ker}~\mathrm{div}_A$ are unitary and
$X_A=U_\pm\mathrm{curl}U_\pm^{-1}$, i.e., the operators $X_A$ and $X_0=\mathrm{curl}$ are
unitary equivalent
\end{theorem}

\pr
The key observation in the proof is that the square of the operator $S_A$ is
diagonal in the sense that $S_A^2:\mathcal{H}^0\rightarrow \mathcal{H}^0$ and
$S_A^2:\mathcal{H}^1\rightarrow \mathcal{H}^1$. Indeed, from formula (\ref{sa}) we
immediately have
\begin{equation}\label{sa2}
S_A^2=\left(
\begin{array}{cc}
P_A{\rm curl}_AP_A{\rm curl}_AP_A-d_A\mathrm{div}_A & 0 \\
 0 & -{\rm div}_Ad_A
\end{array}
\right).
\end{equation}

By the invariance principle for wave operators (see \cite{RS3}) the wave
operators $W_\pm(S_0^2,S_A^2)$ for the pair $S_0^2,S_A^2$ also exist and are
unitary. Moreover, we have the following formula which relates
$W_\pm(S_0,S_A)$ and $W_\pm(S_0^2,S_A^2)$:
\begin{eqnarray}\label{ww2}
\left( W_\pm(S_0,S_A)-W_\pm(S_0^2,S_A^2)\right) E(\mathbb{R}_+)=0, \\
\left( W_\pm(S_0,S_A)-W_\mp(S_0^2,S_A^2)\right) E(\mathbb{R}_-)=0, \nonumber
\end{eqnarray}
where $E$ is the spectral measure of $S_0$.

Since $\mathrm{Ker}~\mathrm{div}$ is an invariant subspace for $S_0$ this subspace is
also invariant for the spectral measure of $S_0$, $E(\mathbb{R}_\pm):
\mathrm{Ker}~\mathrm{div}\rightarrow \mathrm{Ker}~\mathrm{div}$, and, in
particular, $E(\mathbb{R}_\pm): \mathrm{Ker}~\mathrm{div}\rightarrow \mathcal{H}^1$. But the
operator $S_A^2$ is diagonal, and hence $W_\pm(S_0^2,S_A^2):
\mathcal{H}^1\rightarrow \mathcal{H}^1$. The last two observations together
with (\ref{ww2}) show that
\begin{equation}\label{wprop}
W_\pm(S_0,S_A):\mathrm{Ker}~\mathrm{div}\rightarrow
\mathcal{H}^1.
\end{equation}

Now using (\ref{wprop}), the definition of the operator $S_A$ and the intertwining property of
wave operators we have
for any $\omega \in \mathrm{Ker}~\mathrm{div}\bigcap \mathfrak D (S_0)$
\begin{equation}\label{wcalc}
\begin{array}{l}
W_\pm(S_0,S_A)S_0(\omega,0)=S_AW_\pm(S_0,S_A)(\omega,0),\\
\\
W_\pm(S_0,S_A)(P_0\mathrm{curl}~\omega,0)=S_A(W_\pm(S_0,S_A)P_0\omega,0),\\
\\
(W_\pm(S_0,S_A)P_0\mathrm{curl}~\omega,0)=
(X_AW_\pm(S_0,S_A)P_0\omega,\mathrm{div}_AW_\pm(S_0,S_A)P_0\omega).
\end{array}
\end{equation}

From the last line in (\ref{wcalc}) we infer that
$\mathrm{div}_AW_\pm(S_0,S_A)P_0\omega=0$. Therefore
\begin{equation}\label{imw}
\mathrm{Im}~W_\pm(S_0,S_A)P_0\subset
\mathrm{Ker}~\mathrm{div}_A,
\end{equation}
and since $W_\pm(S_0,S_A)$ are unitary
operators we also have
$$
U_\pm^*U_\pm=P_0W_\pm(S_0,S_A)^*W_\pm(S_0,S_A)P_0=P_0.
$$

A calculation for $W_\pm(S_0,S_A)^*$
similar to (\ref{wcalc}) shows that $U_\pm U_\pm^*=P_A$. Therefore $U_\pm$
are unitary operators.

Finally (\ref{imw}) and the last equality in (\ref{wcalc}) imply that
$$
P_AW_\pm(S_0,S_A)P_0\mathrm{curl}=
X_AP_AW_\pm(S_0,S_A)P_0,
$$
or $U_\pm\mathrm{curl}=X_AU_\pm$. This completes the proof of the theorem.

\qed

\section{An extension of the main result}

In conclusion we make a few remarks on possible extensions of Theorems \ref{sauniteq}
and \ref{xauniteq} for connections $A$ which are not smooth.

First, instead of the space $\Omega^1(\mathbb{R}^{3},\k)$ of smooth connection one-forms we shall
consider the space
 ${\mathcal W}^1_{2,loc}({\mathbb{R}}^3,\mathfrak k)$ of $\k$-valued one-forms on
$\mathbb{R}^3$ whose coefficients are elements of the Sobolev space
${W}^1_{2,loc}({\mathbb{R}}^3)$.

One can also define the gauge group ${\mathcal W}^2_{2,loc}({\mathbb{R}}^3,K)$ of Sobolev gauge
transformations by requiring that $g\in {\mathcal W}^2_{2,loc}({\mathbb{R}}^3,K)$ if and only
if $dgg^{-1}\in {\mathcal W}^1_{2,loc}({\mathbb{R}}^3,\mathfrak k)$ (see, for instance, \cite{JT}).
The space ${\mathcal W}^2_{2,loc}({\mathbb{R}}^3,K)$ is indeed a group continuously acting
by gauge transformations on ${\mathcal W}^1_{2,loc}({\mathbb{R}}^3,\mathfrak k)$.

Obviously ${\mathcal W}^2_{2,loc}({\mathbb{R}}^3,K)$ also acts by unitary transformations on the
space $\mathcal{H}^1$.

Now from Theorems \ref{sauniteq} and \ref{xauniteq}  we infer the following statement.

\begin{theorem}
Let $A\in {\mathcal W}^1_{2,loc}({\mathbb{R}}^3,\mathfrak k)$ ,
$F=dA+\frac{1}{2}[A\wedge A]$ the curvature of $A$, where the derivatives of the coefficients of
$A$ should be understood in the sense of generalized functions.
There exists a positive constant $C$ such that if for some $s>3$
and $g\in {\mathcal W}^2_{2,loc}({\mathbb{R}}^3,K)$
\begin{eqnarray*}
\|g\circ F(1+|x|^2)^{s}\|_{(\infty)}< C,  \\
\|g\circ A(1+|x|^2)^{\frac{s}{2}}\|_{(\infty)}< C,
\end{eqnarray*}
then

(i) The corresponding operator $S_A$ is selfadjoint with the domain $g^{-1}\mathfrak D (S_0)$.

(ii) The wave operators
$$
W_{\pm}(S_0,S_A)=s-\lim_{t\rightarrow \pm \infty}e^{iS_At}e^{-iS_0t}
$$
for the pair $S_0,S_A$ exist and are unitary. In particular,
the selfadjoint operators $S_0$ and $S_A$ are unitary equivalent.

(iii) The operators $U_\pm=P_AW_\pm(S_0,S_A)P_0:\mathrm{Ker}~\mathrm{div}\rightarrow
\mathrm{Ker}~\mathrm{div}_A$ are unitary and for the operator $X_A=P_AX_AP_A$ we have
$X_A=U_\pm\mathrm{curl}U_\pm^{-1}$, i.e., the operators $X_A$ and $X_0$ are
unitary equivalent
\end{theorem}


\begin{thebibliography}{99}

\bibitem{Ag} Agmon, S., Spectral properties of Schr\"{o}dinger
operators and scattering theory, {\em Ann. Scuola Norm. Sup. Pisa Cl. Sci.},
{\bf (4)  2}  (1975), no. 2, 151--218.

\bibitem{BS1} Birman, M. Sh., Solomyak, M. Z., $L_2$-theory of the Maxwell
operator in arbitrary domains, {\em Russian Math. Surveys}, {\bf 42} (1987),
no. 6, 75--96

\bibitem{BS2} Birman, M. Sh., Solomyak, M. Z., The self-adjoint Maxwell
operator in arbitrary domains, {\em Algebra and Analysis}, {\bf 1} (1989),
no. 1, 96--110.

\bibitem{Fi} Filonov, N., Spectral analysis of the selfadjoint operator curl
in a region of finite measure, {\em St. Petersburg Math. J.}, {\bf 11}
(2000), 1085--1095.

\bibitem{JT} Jaffe, A., Taubes, C. H., Vortices and monopoles,
Birkh\"{a}user, Boston (1980).

\bibitem{JK} Jensen, A., Kato, T., Spectral properties of Schr\"{o}dinger
operators and time--decay of the wave functions, {\em Duke Math. J.}, {\bf
46} (1979), 583--611.

\bibitem{Kato1} Kato, T., Wave operators and similarity for some
non-selfadjoint operators, {\em Math. Annalen}, {\bf 162} (1966), 258--279.

\bibitem{LU} Ladyzenskaya, O. A., Uraltseva, N. N., Linear and quasilinear
equations of elliptic type, Nauka, Moscow (1973).

\bibitem{LS} Leinfelder, H., Simander, C., Schr\"{o}dinger operators with
singular magnetic vector potentials, {\em Math. Z.}, {\bf 176} (1981),
1--19.

\bibitem{P}Picard, R., On a selfadjoint realization of curl in exterior domains,
{\em Math. Z.}  {\bf 229}  (1998), 319--338.

\bibitem{RS2} Reed, M., Simon, B., Methods of modern mathematical physics,
vol. 2, Academic Press, London (1975).

\bibitem{RS3} Reed, M., Simon, B., Methods of modern mathematical physics,
vol. 3, Academic Press, London (1979).

\bibitem{RS4} Reed, M., Simon, B., Methods of modern mathematical physics,
vol. 4, Academic Press, London (1978).


\end{thebibliography}
\end{document}